\documentclass[a4paper,british]{amsart}
\usepackage{babel}
\usepackage{amssymb}

\newcommand{\heute}{5 June 2002}

\numberwithin{equation}{section}

\newtheorem{theorem}[equation]{Theorem}
\newtheorem{lemma}[equation]{Lemma}
\newtheorem{corollary}[equation]{Corollary}

\theoremstyle{definition}

\newtheorem{definition}[equation]{Definition}
\newtheorem{conjecture}[equation]{Conjecture}

\theoremstyle{remark}
\newtheorem{remark}[equation]{Remark}

\newcommand{\f}[1][p]{\mathbb{F}_{#1}}

\newcommand{\Bild}{\operatorname{Im}}
\newcommand{\Inf}{\operatorname{Inf}}

\newcommand{\Res}{\operatorname{Res}}

\newcommand{\Tr}{\operatorname{Tr}}
\newcommand{\Ann}{\operatorname{Ann}}

\newcommand{\GL}{\text{\sl GL}}
\newcommand{\GLn}[1][p]{\GL(\f)}
\newcommand{\A}{\mathcal{A}}
\newcommand{\calH}{\mathcal{H}}
\newcommand{\calK}{\mathcal{K}}
\newcommand{\calS}{\mathcal{S}}

\newcommand{\N}{\mathcal{N}}

\newcommand{\coho}[2][*]{\mathrm{H}^{#1}(#2)}
\newcommand{\PH}[1]{\mathrm{PH}^*(#1)}

\newcommand{\rp}{\mathcal{P}}

\newcommand{\frakp}{\mathfrak{p}}
\newcommand{\frakq}{\mathfrak{q}}

\newcommand{\tauH}[1][\calH]{\tau_{#1}}
\newcommand{\saH}[1][a]{s_{\mathfrak{#1},\calH}}
\newcommand{\tauaH}[1][a]{\tau_{\mathfrak{#1},\calH}}
\newcommand{\taua}[1][a]{\tau_{\mathfrak{#1}}}
\newcommand{\sa}[1][a]{s_{\mathfrak{#1}}}
\newcommand{\SaH}[1][a]{\calS_{\mathfrak{#1}, \calH}}
\newcommand{\fui}[1][i]{f_{\underline{#1}}}
\newcommand{\fun}{\fui[n]}
\newcommand{\fuo}{\fui[1]}

\begin{document}

\title{On Carlson's depth conjecture in group cohomology}

\author[D.~J. Green]{David J. Green}
\address{Dept of Mathematics \\ Univ.\@ of Wuppertal \\
D--42097 Wuppertal \\ Germany}
\email{green@math.uni-wuppertal.de}

\date{\heute}

\begin{abstract}
We establish a weak form of Carlson's conjecture on the depth of the
mod-$p$ cohomology ring of a $p$-group. In particular,
Duflot's lower bound for the depth is tight if and only if the cohomology
ring is not detected on a certain family of subgroups.
The proofs use the structure of the cohomology ring as a comodule over the
cohomology of the centre via the multiplication map. We demonstrate the
existence of systems of parameters (so-called polarised systems)
which are particularly well adapted to this comodule structure.
\end{abstract}

\maketitle

\section*{Introduction}
\noindent
Let $G$ be a finite $p$-group and $C$ its greatest central elementary
abelian subgroup. Write $k$ for the prime field $\f$. Cohomology will
always be with coefficients in~$k$.
Denote by $r$ the $p$-rank of~$G$, and by $z$ the $p$-rank of $C$.
Following Broto and Henn~\cite{BrHe:Central}
we shall exploit the fact that the multiplication map
$\mu \colon G \times C \rightarrow G$, $(g,c) \mapsto g.c$ is a group
homomorphism. The main result of this paper is as follows:

\begin{theorem}
\label{theorem:introED1}
Suppose that $G$ is a $p$-group whose centre has $p$-rank $z$.
Then the following statements are equivalent:
\begin{enumerate}
\item
The mod-$p$ cohomology ring $\coho{G}$ is not detected on the
centralizers of its rank $z+1$ elementary abelian subgroups.
\item
There is an associated prime $\frakp$ such that
$\coho{G}/\frakp$ has dimension $z$.
\item
The depth of $\coho{G}$ equals $z$.
\end{enumerate}
\end{theorem}

\noindent
This is a special case of a conjecture due to
Carlson~\cite{Carlson:DepthTransfer},
reproduced here as Conjecture~\ref{conjecture:ED}. Recall that Duflot proved
in~\cite{Duflot:Depth} that $z$ is a lower bound for the depth. So this
result characterises the cases where Duflot's lower bound is tight.

Theorem~\ref{theorem:introED1} is proved as Theorem~\ref{theorem:ED1}.
The proof rests upon the existence of \emph{polarised systems}:
homogeneous systems of parameters for $\coho{G}$ which are particularly well
adapted to the $\coho{C}$-comodule structure.
There are two extreme types of behaviour
which a cohomology class $x \in \coho{G}$ can demonstrate under
the comodule structure map $\mu^*$: one extreme is that
$\Res_C(x)$ is nonzero, and so $\mu^*(x) = 1 \otimes \Res_C(x) + \text{other
terms}$. The other extreme is that $x$ is primitive, meaning that
$\mu^*(x) = x \otimes 1$. Roughly speaking, a polarised system of parameters
is one consisting solely of elements which each exhibit one or the other
of these extreme kinds of behaviour. The precise definition,
which ensures that each such system is a detecting sequence for the
depth of $\coho{G}$, is given in Definition~\ref{definition:polarised}\@.
Polarised systems of parameters always exist, as is proved in
Theorem~\ref{theorem:existence}\@.
\medskip

\noindent
In addition to Theorem~\ref{theorem:introED1} we also prove a weak form
of the general case of Carlson's conjecture.
This is done in Theorem~\ref{theorem:polarisedEqualities},
which includes the following statement:

\begin{theorem}
\label{theorem:introMyEDgen}
Let $\zeta_1$, \dots, $\zeta_z$, $\kappa_1$, \dots, $\kappa_{r-z}$ be
a polarised system of parameters for $\coho{G}$.
Then $\coho{G}$ has depth $z + \sa$, where
$\sa \in \{0,\ldots,r-z\}$
is defined by
\[
\sa = \max \{ s \leq r-z \mid \text{$\kappa_1,\ldots,\kappa_s$
is a regular sequence in $\coho{G}$} \} \, .
\]
\end{theorem}

\section{Primitive comodule elements}
\label{section:primitive}
\noindent
Group multiplication $\mu \colon G \times C \rightarrow G$ is a
group homomorphism. As observed by Broto and Henn~\cite{BrHe:Central},
this means that $\coho{G}$ inherits the structure of a comodule
over the coalgebra $\coho{C}$.

Recall that $x \in \coho{G}$ is called a primitive comodule element
if $\mu^*(x) = x \otimes 1 \in \coho{G \times C} \cong \coho{G} \otimes_k
\coho{C}$. As the comodule structure map $\mu^*$ is simultaneously
a ring homomorphism, it follows that the primitives form a subalgebra
$\PH{G}$ of $\coho{G}$. As the quotient map $G \rightarrow G/C$ coequalises
$\mu$ and projection onto the first factor of $G \times C$, it follows
that the image of inflation from $\coho{G/C}$ is contained in
$\PH{G}$.

\begin{lemma}
\label{lemma:quotientComodule}
Suppose that $I$ is a homogeneous ideal in $\coho{G}$ which is generated
by primitive elements. Then
\[
\mu^*(I) \subseteq I \otimes_k \coho{C} \, ,
\]
and so $\mu^*$ induces a ring homomorphism
\[
\lambda \colon \coho{G}/I \rightarrow \coho{G}/I \otimes_k \coho{C}
\]
which induces an $\coho{C}$-comodule structure on $\coho{G}/I$.
\end{lemma}

\begin{proof}
For $x \in I$ and $y \in \coho{G}$ one has
$\mu^*(xy) = \mu^*(x) \mu^*(y) = (x \otimes 1) \mu^*(y)
\in I \otimes_k \coho{C}$.
\end{proof}

\begin{lemma}
\label{lemma:myBrotoCarlsonHenn}
Suppose that $\zeta_1, \ldots, \zeta_t$ is a sequence of homogeneous
elements of $\coho{G}$ whose restrictions form a regular sequence
in $\coho{C}$, and suppose that $I$~is an ideal in $\coho{G}$ generated
by primitive elements. Then $\zeta_1, \ldots, \zeta_t$
is a regular sequence for the quotient ring $\coho{G}/I$.
\end{lemma}

\begin{proof}
Carlson's proof for the case $I=0$ (Proposition 5.2 of~\cite{Carlson:Problems})
generalises easily. Denote by $R$ the polynomial algebra
$k[\zeta_1, \ldots, \zeta_t]$.
The map $\lambda$~of Lemma~\ref{lemma:quotientComodule}
induces an $R$-module structure on $\coho{G}/I \otimes_k \coho{C}$,
and $\lambda$~is a split monomorphism of $R$-modules, the splitting map
being induced by projection onto the first factor $G \times C \rightarrow G$.
So as an $R$-module $\coho{G}/I$ is a direct summand
of $\coho{G}/I \otimes_k \coho{C}$. The result will therefore follow
if we can show that
$\coho{G}/I \otimes_k \coho{C}$ is a free $R$-module.

To see that $\coho{G}/I \otimes_k \coho{C}$ is indeed a free $R$-module
set
\[
F_i := \sum_{j \geq i} \coho[j]{G}/(I \cap \coho[j]{G})
\]
and observe that $F_i \otimes_k \coho{C}$ is an $R$-submodule of
$\coho{G}/I \otimes_k \coho{C} = F_0 \otimes_k \coho{C}$.
Projection $G \times C \rightarrow C$ makes
$F_i \otimes_k \coho{C} / F_{i+1} \otimes_k \coho{C}$ a free
$\coho{C}$-module.
Now for $x \in F_i$, $y \in \coho{C}$ and $\theta \in R$ we have
$\theta.(x \otimes y) \in x \otimes (\Res^G_C \theta). y
+ F_{i+1} \otimes_k \coho{C}$. 
So the $R$-module structure on
$F_i \otimes_k \coho{C} / F_{i+1} \otimes_k \coho{C}$ is induced by
the restriction map $R \rightarrow \coho{G} \rightarrow \coho{G}/I
\rightarrow \coho{C}$ from the free $\coho{C}$-structure.
But $\coho{C}$ is a free $R$-module by Theorem 10.3.4 of \cite{Evens:book},
because the restrictions of $\zeta_1,\ldots,\zeta_t$ form a regular
sequence. So
$F_i \otimes_k \coho{C} / F_{i+1} \otimes_k \coho{C}$ is a free $R$-module
for all~$i$, whence it follows that
$F_0 \otimes_k \coho{C} / F_i \otimes_k \coho{C}$ is a free $R$-module for
all~$i$. As the degree of each homogeneous element of $F_i \otimes_k \coho{C}$
is at least~$i$, it follows that $\coho{G}/I \otimes_k \coho{C}$ is
itself a free $R$-module.
\end{proof}

\begin{corollary}
\label{coroll:myBrotoCarlsonHenn}
Let $G$ be a $p$-group whose centre has $p$-rank~$z$.
Suppose that there is a length~$s$ regular sequence in $\coho{G}$
which consists entirely of primitive elements. Then the depth
of $\coho{G}$ is at least $z + s$.
\end{corollary}

\begin{proof}
Let $I$ be the ideal generated by the primitive elements in the
regular sequence. Let $\zeta_1,\ldots,\zeta_z$ be elements of
$\coho{G}$ whose restrictions form a homogeneous system of
parameters for $\coho{C}$: it is well known that such sequences exist.
Now apply Lemma~\ref{lemma:myBrotoCarlsonHenn}.

One concrete example of such classes $\zeta_i$ is obtained as follows.
Let $\rho_G$ be the regular representation of~$G$, and $\zeta_i$
the Chern class $c_{p^n - p^{n-i}}(\rho_G)$ for $1 \leq i \leq z$,
where $p^n$ is the order of~$G$. Then $\rho_G$ restricts to~$C$
as $|G:C|$ copies of the regular representation $\rho_C$,
whence $\Res_C(\zeta_i) = c_{p^z - p^{z-i}}(\rho_C)^{|G:C|}$.
But the $c_{p^z - p^{z-i}}(\rho_C)$ are (up to a sign)
the Dickson invariants.
See the proof of Theorem~\ref{theorem:existence} for more details.
\end{proof}

\section{Polarised systems of parameters}
\label{section:polarised}
\noindent
We shall now give the definition of a polarised system of parameters, the
key definition of this paper. In fact we shall introduce two closely related
concepts: the axioms for a polarised system
(Definition~\ref{definition:polarised}) are easily checked in practice,
whereas the special polarised systems of
Definition~\ref{definition:specialPolarised} have precisely the properties
we shall need to investigate depth. Lemma~\ref{lemma:polarisedDefinitions}
shows that the two concepts are more or less interchangeable.

\begin{definition}
\label{definition:ACG}
Let $G$ be a $p$-group of $p$-rank $r$ whose centre has $p$-rank $z$.
Denote by $C$ the greatest central elementary abelian subgroup of~$G$, and
set
\begin{align*}
\A^C(G) & :=
\{ V \leq G \mid \text{$V$ is elementary abelian and contains~$C$} \} \, , \\
\A^C_d(G) & :=
\{ V \in \A^C(G) \mid \text{$V$ has $p$-rank $d$} \} \, \\
\calH^C_d(G) & := \{ C_G(V) \mid V \in \A^C_d(G) \} \, .
\end{align*}
So $\A^C(G)$ is the disjoint union of the $\A^C_{z+s}(G)$ for
$0 \leq s \leq r-z$.
\end{definition}

\begin{definition}
\label{definition:polarised}
Let $G$ be a $p$-group of $p$-rank $r$ whose centre has $p$-rank $z$.
Recall that inflation map $\Inf \colon \coho{V/C} \rightarrow \coho{V}$
is a split monomorphism for each $V \in \A^C(G)$,
and so its image $\Bild \Inf$ is isomorphic to $\coho{V/C}$.

A system $\zeta_1,\ldots,\zeta_z$, $\kappa_1,\ldots,\kappa_{r-z}$
of homogeneous elements of $\coho{G}$ shall be called a polarised
system of parameters if the following four axioms are satisfied.
\begin{description}
\item[(PS1)]
$\Res_C(\zeta_1)$, \dots, $\Res_C(\zeta_z)$ is a system of parameters for
$\coho{C}$.
\item[(PS2)]
$\Res_V(\kappa_j)$ lies in $\Bild \Inf$ for each $1 \leq j \leq r-z$ and for
each $V \in \A^C(G)$.
\item[(PS3)]
For each $V \in \A^C(G)$,
the restrictions $\Res_V(\kappa_1), \ldots, \Res_V(\kappa_s)$
constitute a system of
parameters for $\Bild \Inf$. Here $z+s$ is the rank of~$V$.
\item[(PS4)]
$\Res_V(\kappa_j) = 0$ for $V \in \A^C_{z+s}(G)$ with $0 \leq s < j \leq r-z$.
\end{description}
\end{definition}

\begin{remark}
Polarised systems of parameters always exist, as we shall see in
Theorem~\ref{theorem:existence}\@.
Observe that Axiom (PS1) involves only the $\zeta_i$, whereas the remaining
axioms involve only the $\kappa_j$.
Basically Axiom (PS1) says that $\zeta_1,\ldots,\zeta_z$ is a regular
sequence which can be detected on the centre, (PS2) says that the $\kappa_j$
are primitive after raising to suitably high $p$th powers, and (PS3) says
that the $\kappa_j$ together with the $\zeta_i$ will form a detecting
sequence for the depth of $\coho{G}$. Axiom (PS4) is a more technical
condition which we shall only use once: it is needed in
Lemma~\ref{lemma:polarisedDefinitions} to show that, after raising to
a suitably high $p$th power, each $\kappa_j$ is a sum of transfer classes
as required by Axiom (PS5) below.
\end{remark}

\begin{lemma}
Polarised systems of parameters for $\coho{G}$ are indeed systems of parameters.
\end{lemma}

\begin{proof}
Let $V \in \A^C_{z+s}(G)$.
The restrictions of $\zeta_1,\ldots,\zeta_z, \kappa_1,\ldots, \kappa_s$
constitute a system of parameters for~$\coho{V}$ by (PS1) and (PS3)\@.
Hence $\zeta_1,\ldots,\zeta_z, \kappa_1, \ldots, \kappa_{r-z}$ are
algebraically independent over~$k$, for we may choose $V$ to have $p$-rank~$r$.

Now let $\gamma$ be a homogeneous element of $\coho{G}$. For $V \in \A^C(G)$
there is a monic polynomial $f_V(x)$
with coefficients in $k[\zeta_1,\ldots,\zeta_z, \kappa_1, \ldots, \kappa_{r-z}]$
such that $f_V(\gamma)$ has zero restriction to~$V$. Taking the product
of all such polynomials one obtains a polynomial $f(x)$ such that
$f(\gamma)$ has zero restriction to each maximal elementary abelian subgroup
of~$G$. So $f(\gamma)$ is nilpotent by a well-known result of Quillen.
\end{proof}

\begin{definition}
\label{definition:specialPolarised}
Let $G$ be a $p$-group of $p$-rank $r$ whose centre has $p$-rank $z$.
A system $\zeta_1,\ldots,\zeta_z$, $\kappa_1,\ldots,\kappa_{r-z}$
of homogeneous elements of $\coho{G}$ shall be called a special polarised
system of parameters if it satisfies the following five axioms:
(PS1), (PS3), (PS4) and
\begin{description}
\item[(PS$\mathbf{2'}$)]
$\kappa_j$ is a primitive element of the $\coho{C}$-comodule $\coho{G}$
for each $1 \leq j \leq r-z$.
\item[(PS5)]
$\kappa_j$ lies in
$\sum_{H \in \calH^C_{z+i}(G)} \Tr^G_H(\coho{H})$ for each
$1 \leq i \leq j \leq r-z$.
\end{description}
\end{definition}

\begin{lemma}
\label{lemma:polarisedDefinitions}
Axiom (PS$2'$) implies Axiom (PS2), and so every special polarised
system is a polarised system of parameters.
Conversely for each polarised system
$\zeta_1,\ldots,\zeta_z$, $\kappa_1, \ldots, \kappa_{r-z}$
there is a nonnegative integer $N$ such that $\zeta_1,\ldots,\zeta_z$,
$\kappa_1^{p^N}, \ldots, \kappa_{r-z}^{p^N}$ is a special polarised
system of parameters for $\coho{G}$.
\end{lemma}

\begin{proof}
Let $V \in \A^C(G)$.
Restriction from $\coho{G}$~to $\coho{V}$ is a map of $\coho{C}$-comodules
and so sends primitive elements to primitive elements. But the subalgebra
of primitive elements of $\coho{V}$ coincides with the image
of inflation from $V/C$.  So (PS$2'$) implies (PS2).

Now suppose $\zeta_1,\ldots,\zeta_z$, $\kappa_1, \ldots, \kappa_{r-z}$ is
a polarised system for $\coho{G}$. For each $1 \leq j \leq r-z$
the restriction of $\kappa_j$ to each $V \in \A^C(G)$ is primitive
by (PS2)\@. Hence the element $\mu^*(\kappa_j) - \kappa_j \otimes 1$
of $\coho{G} \otimes \coho{C}$ has zero restriction to every maximal
elementary abelian subgroup and is therefore nilpotent.

For fixed $1 \leq i \leq r-z$, denote
by $\calK$ the set consisting of those subgroups $K$~of $G$ such
that $C_G(K)$ is not $G$-conjugate to any subgroup of
any $H \in \calH^C_{z+i}$. Following Carlson (Proof of Corollary 2.2~of
\cite{Carlson:DepthTransfer}), observe that every $K \in \calK$ has
$p$-rank less than $z+i$. Moreover every $K \in \calK$ is contained
in $K_C = \langle K,C \rangle$, and $K_C$~itself lies in~$\calK$.
So $\Res_{K_C}(\kappa_j) = 0$ for all $j \geq i$ by (PS4)\@.
Hence each such $\kappa_j$ lies in the
radical ideal $\sqrt{J'}$, where $J'$ is the ideal
$\bigcap_{K \in \calK} \ker \Res_K$. So by Benson's result on the
image of the transfer map (Theorem~1.1 of \cite{Benson:ImTr}) the
$\kappa_j$ also lie in $\sqrt{J}$, where $J$ is the ideal
$\sum_{H \in \calH^C_{z+i}(G)} \Tr^G_H(\coho{H})$.
\end{proof}

\begin{remark}
The above proof is the only time we shall make use of the Axiom (PS4)\@.
In particular, the results of \S\ref{section:specialPolarisedDepth}
do not depend on (PS4)\@.
I do not know whether or not (PS4) is a consequence of the remaining
axioms for a special polarised system of parameters.

Axiom (PS5) will be used in Lemma~\ref{lemma:kappaAssocPrime} to
prove the existence of an associated prime ideal with desirable properties.
\end{remark}

\begin{lemma}
\label{lemma:genBrotoHenn}
Suppose that $\zeta_1,\ldots,\zeta_z$, $\kappa_1,\ldots,\kappa_{r-z}$ is
a polarised system of parameters for $\coho{G}$.
Let $0 \leq s \leq r-z$.
Then the sequence
$\zeta_1, \ldots, \zeta_z, \kappa_1, \ldots, \kappa_s$
is regular in $\coho{G}$ if and only if the sequence
$\kappa_1, \ldots, \kappa_s$ is regular in $\coho{G}$.
\end{lemma}

\begin{proof}
Recall that regular sequences may be permuted at will. Moreover, replacing
one element of a sequence by its $p$th power has no effect on whether the
sequence is regular or not. Hence by Lemma~\ref{lemma:polarisedDefinitions}
it suffices to prove the result for special polarised systems.

So we may assume that the given sequence is a special polarised system
of parameters.
Let $I$ be the homogeneous ideal $I$ in $\coho{G}$ generated
by $\kappa_1, \ldots, \kappa_s$.
By Axiom (PS$2'$) this ideal is generated by primitive elements.
Also, the restrictions of $\zeta_1,\ldots,\zeta_z$ form a regular
sequence in $\coho{C}$ by Axiom~(PS1)\@. Therefore
Lemma~\ref{lemma:myBrotoCarlsonHenn} tells us that $\zeta_1,\ldots,\zeta_z$
constitute a regular sequence for $\coho{G}/I$. 
\end{proof}

\section{Three depth-related numbers}
\label{section:threeNumbers}
\noindent
In this section we shall introduce three numbers $\tauH$, $\taua$~and
$\tauaH$, each of which is an approximation to the depth $\tau$ of $\coho{G}$.

\begin{definition}
Let $G$ be a $p$-group of $p$-rank $r$ whose centre has $p$-rank $z$.
Write $\tau$ for the depth of $\coho{G}$ and set
\[ \tauH := \max \{ d \in \{z, \ldots, r\} \mid
\text{The family $\calH^C_d(G)$ detects $\coho{G}$} \} \, .
\]
\end{definition}

\noindent
In~\cite{Carlson:DepthTransfer} Carlson formulates the following conjecture:
\begin{conjecture}[Carlson]
\label{conjecture:ED}
The number~$\tauH$ coincides with the depth $\tau$~of $\coho{G}$.
Moreover, $\coho{G}$ has an associated prime
$\frakp$ such that $\dim \coho{G}/\frakp = \tau$.
\end{conjecture}

\noindent
In fact, Carlson formulates the conjecture not just for $p$-groups,
but for arbitrary finite groups. In this article however we only
consider $p$-groups.
In Theorem~\ref{theorem:ED1} we shall prove a special case of this
conjecture, after deriving a partial result for the general case in
Theorem~\ref{theorem:polarisedEqualities}. For this we need two
more depth-related numbers.

\begin{definition}
Let $G$ be a $p$-group of $p$-rank $r$ whose centre has $p$-rank $z$, and
let $\mathfrak{a} =
(\zeta_1,\ldots,\zeta_z, \kappa_1,\ldots,\kappa_{r-z})$
be a polarised system of parameters for $\coho{G}$\@.

Define $\taua$ to be $z + \sa$, where $\sa$ is the largest
$s \in \{0, \ldots, r-z\}$ such that $\kappa_1,\ldots, \kappa_s$
is a regular sequence in $\coho{G}$.

Let $\SaH$~be the subset of $\{0, \ldots, r-z\}$ such that $s$~lies
in $\SaH$ if and only if the restriction map
\[
\coho{G}/(\kappa_1,\ldots,\kappa_{s-1}) \rightarrow
\prod_{H \in \calH^C_{z+s}}
\coho{H}/(\Res_H \kappa_1,\ldots, \Res_H \kappa_{s-1})
\]
is injective\footnote{So $\SaH$ always contains~$0$, and $1$~lies in $\SaH$
if and only if the family $\calH^C_{z+1}$ detects $\coho{G}$.}.
Define $\saH := \max \SaH$ and $\tauaH := z + \saH$.
\end{definition}

\begin{lemma}
\label{lemma:tauCdash}
Let $\mathfrak{a} = (\zeta_1,\ldots,\zeta_z$, $\kappa_1,\ldots,\kappa_{r-z})$
be a polarised system of parameters for $\coho{G}$\@.
If $s > 0$ lies in $\SaH$ then
the family $\calH^C_{z + s}(G)$ detects $\coho{G}$
and $s-1$ lies in~$\SaH$.
Therefore $\tauH \geq \tauaH$ and $\SaH = \{ 0, \ldots, \saH \}$.
\end{lemma}

\noindent
For the proof we shall need an elementary fact about regular sequences.

\begin{lemma}
\label{lemma:genGrComm}
Suppose that $R,S$ are connected graded commutative $k$-algebras and
that $f \colon R \rightarrow S$ is an algebra homomorphism which respects
the grading. Suppose further that $\zeta_1,\ldots,\zeta_d$ is a family
of homogenous positive-degree elements of $R$ satisfying the following
conditions:
\begin{enumerate}
\item
$f(\zeta_1),\ldots,f(\zeta_d)$ is a regular sequence in~$S$.
\item
The induced map $f_d \colon R/(\zeta_1,\ldots,\zeta_d)
\rightarrow S/(f(\zeta_1),\ldots,f(\zeta_d))$ is an injection.
\end{enumerate}
Then $f \colon R \rightarrow S$ is an injection.
\end{lemma}

\begin{proof}
It suffices to prove the case $r=1$. Write $\zeta$~for $\zeta_1$.
Let $a \not = 0$ be an element of $\ker(f)$ whose degree is as
small as possible. Since $f(a) = 0$ in $S/(f(\zeta))$ it follows that
there is an $a' \in R$ with $a = a' \zeta$. Since $f(a) = 0$ and $f(\zeta)$
is regular it follows that $a' \in \ker(f)$, contradicting the minimality
of $\deg(a)$.
\end{proof}

\begin{proof}[Proof of Lemma~\ref{lemma:tauCdash}]
Apply Lemma~\ref{lemma:genGrComm} to the family $\kappa_1,\ldots,
\kappa_{s-1}$ with $R = \coho{G}$, $S = \prod_{\calH^C_{z + s}} \coho{H}$
and $f$ the product of the restriction maps.
Because $s \in \SaH$ the induced
map of quotients is an injection. By Axiom~(PS3) the restrictions of
$\kappa_1,\ldots, \kappa_{s-1}$ form a regular sequence in $\coho{V}$
for each $V \in \A^C_{z + s}(G)$, and so by
\cite[Prop.~5.2]{Carlson:Problems} they form a regular sequence in $\coho{H}$
for each $H \in \calH^C_{z + s}$.
Hence the restrictions form a regular sequence in~$S$, and so
the family $\calH^C_{z + s}$ detects $\coho{G}$.
If instead we just invoke the first step in the proof of
Lemma~\ref{lemma:genGrComm}, we see that
the $\coho{H}/(\Res \kappa_1,\ldots,\Res \kappa_{s-2})$ with
$H \in \calH^C_{z + s}$ detect $\coho{G}/(\kappa_1,\ldots,\kappa_{s-2})$.
\end{proof}

\section{Depth and special polarised systems}
\label{section:specialPolarisedDepth}
\noindent
The following fact from commutative algebra is well known.

\begin{lemma}
\label{lemma:assocPrime}
Let $A$ be a graded commutative ring and $M$ a Noetherian graded $A$-module.
Suppose that $\frakp$ is an associated prime of~$M$, and that
$\zeta_1,\ldots,\zeta_d$ is a regular sequence for $M$.
Then $M/(\zeta_1,\ldots,\zeta_d)M$ has an associated prime~$\frakq$ containing
$\frakp$.
\end{lemma}

\begin{proof}
It suffices to prove the case $d=1$. Write $\zeta$~for $\zeta_1$.
Pick $x \in M$ with $\Ann_A(x) = \frakp$. If $x$ lies in $\zeta M$ then
there is an $x' \in M$ with $\zeta x' = x$. As $\zeta$ is regular
it follows that $\Ann_A(x') = \frakp$ too, so replace $x$~by $x'$.
This can only happen finitely often, as $M$ is Noetherian and
$Ax$ is strictly contained in $Ax'$. So we may assume that $x$~does not
lie in $\zeta M$, which means that the image of $x$~in $M/\zeta M$ is
nonzero and annihilated by $\frakp$.
\end{proof}

\begin{lemma}
\label{lemma:kappaAssocPrime}
Let
$\mathfrak{a} = (\zeta_1,\ldots,\zeta_z, \kappa_1,\ldots,\kappa_{r-z})$
be a special polarised system of parameters
for $\coho{G}$, and $I$ the
ideal generated by
the $\kappa_j$ with $j \leq \saH$.
Then the $\coho{G}$-module
$\coho{G}/I$ has an associated prime $\frakp$
which contains
$\kappa_1$, \dots, $\kappa_{r-z}$.
\end{lemma}

\begin{proof}
Set $s = \saH + 1$. The result is trivial if $\saH = r-z$,
so we may assume that $s \leq r-z$.
By definition of $\saH$, the family $\calH = \calH^C_{z + s}$ does not detect
the quotient $\coho{G}/I$.
Pick a class $x \in \coho{G}$ which does not lie in the ideal
$I$,
but whose restriction to each $H \in \calH$ does lie in the ideal
$\coho{H} . \Res_H(I)$.
Let $A$ be the ideal of classes in $\coho{G}$ which annihilate
the image of $x$~in the quotient
$\coho{G}/I$.

For any $j \geq s$ we have $\kappa_j \in \sum_{H \in \calH}
\Tr^G_H \coho{H}$ by Axiom~(PS5),
say $\kappa_j = \sum_H \Tr_H \gamma_H$.
So $\kappa_j x = \sum_H \Tr_H (\gamma_H \Res_H(x))$
by Frobenius reciprocity. Now by assumption $\Res_H(x)$
lies in the ideal generated by $\Res_H(\kappa_1)$, \dots,
$\Res_H(\kappa_{s-1})$; and this by a second application of Frobenius
reciprocity means that $\kappa_j x$ lies in the
ideal $I$. So
$\kappa_j \in A$ for all $j \geq s$, which means that
the $\coho{G}$-module $\coho{G}/I$ has
an associated prime $\frakp$ containing $\kappa_1,\ldots,\kappa_{r-z}$.
\end{proof}

\begin{theorem}
\label{theorem:specialPolarisedEqualities}
Let $G$ be a $p$-group of $p$-rank $r$ whose centre has $p$-rank $z$, and let
$ \mathfrak{a} = (\zeta_1,\ldots,\zeta_z, \kappa_1,\ldots,\kappa_{r-z})$
be a special polarised system of parameters for $\coho{G}$.
Then the numbers $\taua$~and $\tauaH$ both coincide with the depth
$\tau$~of $\coho{G}$.
\end{theorem}

\begin{proof}
We shall prove that $\tau \geq \taua \geq \tauaH = \tau$.
Each $\kappa_j$ is primitive by Axiom~(PS$2'$), so
$\tau \geq \taua$ by Corollary~\ref{coroll:myBrotoCarlsonHenn}.

$\taua \geq \tauaH$:
Suppose that $s \in \SaH$
and $\kappa_1$, \dots, $\kappa_{s-1}$ is a regular sequence in
$\coho{G}$. If $\kappa_1,\ldots,\kappa_s$ is not a regular sequence
then there is some nonzero
$x \in \coho{G}/(\kappa_1,\ldots,\kappa_{s-1})$ annihilated by
$\kappa_s$.
Since $s \in \SaH$ there is some $H \in \calH^C_{z + s}$
such that $\Res_H(x)$ is nonzero in
$\coho{H}/(\Res_H \kappa_1,\ldots,\Res_H \kappa_{s-1})$. But this
is a contradiction since (as in the proof of Lemma~\ref{lemma:tauCdash})
the restrictions of $\kappa_1,\ldots,\kappa_s$
form a regular sequence in $\coho{H}$. By induction on~$s$ we deduce that
$\taua \geq \tauaH$.

$\tauaH = \tau$:
Set $s = \saH$ and denote by $I$ the ideal
$(\kappa_1,\ldots,\kappa_s)$ of $\coho{G}$.
By Lemma~\ref{lemma:kappaAssocPrime} the $\coho{G}$-module $\coho{G}/I$
has an associated prime~$\frakp$ which contains $\kappa_1$, \dots,
$\kappa_{r-z}$.
As $\taua \geq \tauaH$ the sequence
$\kappa_1,\ldots,\kappa_s$ is regular in $\coho{G}$,
so by Lemma~\ref{lemma:genBrotoHenn} the sequence
$\zeta_1,\ldots,\zeta_z$, $\kappa_1$, \dots,
$\kappa_s$ is regular in $\coho{G}$. Therefore the
sequence $\zeta_1,\ldots,\zeta_z$ is regular for $\coho{G}/I$,
and so (by Lemma~\ref{lemma:assocPrime})
the $\coho{G}$-module
$\coho{G}/(\zeta_1,\ldots,\zeta_z,\kappa_1,\ldots,\kappa_s)$
has an associated prime $\frakq$ containing all elements of the homogeneous
system of parameters $\zeta_1$, \dots, $\zeta_z$, $\kappa_1$, \dots,
$\kappa_{r-z}$ for $\coho{G}$. So the depth of this quotient module is zero.
But every regular sequence in $\coho{G}$ can be
extended to a length~$\tau$ regular sequence
(see \cite[\S4.3--4]{Benson:PolyInvts}, for example). So
$\tau = z + s = \tauaH$.
\end{proof}

\section{Depth and polarised systems}
\label{section:polarisedDepth}
\noindent
In this section we shall remove the requirement in
Theorem~\ref{theorem:specialPolarisedEqualities} that the polarised system
be special.

\begin{theorem}
\label{theorem:polarisedEqualities}
Let $G$ be a $p$-group of $p$-rank $r$ whose centre has $p$-rank $z$, and let
$ \mathfrak{a} = (\zeta_1,\ldots,\zeta_z, \kappa_1,\ldots,\kappa_{r-z})$
be a polarised system of parameters for $\coho{G}$.
Then the numbers $\taua$~and $\tauaH$ both coincide with the depth
$\tau$~of $\coho{G}$.
\end{theorem}

\noindent
For the proof we shall need one further fact about regular sequences.

\begin{lemma}
\label{lemma:liftingInjections}
Suppose that $R,S$ are connected graded commutative $k$-algebras and
that $f \colon R \rightarrow S$ is an algebra homomorphism which respects
the grading. Suppose further that $\zeta_1,\ldots,\zeta_d$ is a family
of homogenous positive-degree elements of $R$ satisfying the following
conditions:
\begin{enumerate}
\item
$\zeta_1,\ldots,\zeta_d$ is a regular sequence in~$R$.
\item
The induced map
$\fun \colon R/(\zeta_1^{n_1},\ldots,\zeta_d^{n_d})
\rightarrow S/(f(\zeta_1)^{n_1},\ldots,f(\zeta_d)^{n_d})$
is an injection for certain positive integers $n_1,\ldots,n_d$.
\end{enumerate}
Then the induced map
$\fuo \colon R/(\zeta_1,\ldots,\zeta_d) \rightarrow
S/(f(\zeta_1), \ldots, f(\zeta_d))$ is an injection.
\end{lemma}

\begin{proof}
%
Pick $x \in R$ with $\fuo(x) = 0$.
Then $\fun(\zeta_1^{n_1-1} \ldots \zeta_d^{n_d-1} x) = 0$ and so
\begin{equation}
\label{eqn:peelOff}
\text{$\zeta_1^{n_1-1} \ldots \zeta_d^{n_d-1} x$ lies in
the ideal $(\zeta_1^{n_1}, \ldots, \zeta_d^{n_d})$ of $R$.}
\end{equation}
Set $\zeta' := \zeta_1^{n_1-1} \ldots \zeta_{d-1}^{n_{d-1}-1}$.
Then there are
$a_1, \ldots, a_d \in R$ such that
$\zeta' \zeta_d^{n_d-1} x = \zeta_1^{n_1} a_1 + \cdots + \zeta_d^{n_d} a_d$,
whence
$\zeta_d^{n_d-1} (\zeta' x - \zeta_d a_d) \in (\zeta_1^{n_1},
\ldots, \zeta_{d-1}^{n_{d-1}})$.
As the sequence $\zeta_1^{n_1}, \ldots, \zeta_{d-1}^{n_{d-1}}, \zeta_d^s$ is
regular in~$R$ for $s \geq 1$ we deduce that
$\zeta' x \in (\zeta_1^{n_1}, \ldots, \zeta_{d-1}^{n_{d-1}}, \zeta_d)$.
So we have reduced Eqn.~\eqref{eqn:peelOff} to the case $n_d = 1$
without altering the remaining $n_t$. As regular sequences may be permuted
at will we deduce that $x \in (\zeta_1, \ldots, \zeta_d)$.
\end{proof}

\begin{proof}[Proof of Theorem~\ref{theorem:polarisedEqualities}]
Arguing exactly as in the proof of
Theorem~\ref{theorem:specialPolarisedEqualities} one shows that
$\taua \geq \tauaH$.
By Lemma~\ref{lemma:polarisedDefinitions} there is an integer $N \geq 0$
such that the system of parameters
$\mathfrak{b} = (\zeta_1,\ldots,\zeta_z,
\kappa_1^{p^N}, \ldots, \kappa_{r-z}^{p^N})$ is special polarised.
From the definition one sees that $\taua = \taua[b]$.
So by Theorem~\ref{theorem:specialPolarisedEqualities} it suffices to prove
that $\tauaH \geq \tauaH[b]$. But this is a consequence of
Lemma~\ref{lemma:liftingInjections} applied to the sequence
$\kappa_1, \ldots, \kappa_{\saH[b]-1}$ with $R = \coho{G}$,
$S = \prod_{\calH} \coho{H}$ and $f$~the restriction map, where
$\calH = \calH^C_{z+\saH[b]}(G)$ and $n_i = p^N$ for all~$i$.
\end{proof}

\section{Dickson invariants}
\label{section:Dickson}
\noindent
Let $V$ be an $m$-dimensional $k$-vector space. We shall make extensive
use of the Dickson invariants, the polynomial generators of the ring
of $\GL(V)$-invariants in $k[V]$. See Benson's
book~\cite{Benson:PolyInvts} for proofs of the properties of these invariants.

Denote by $f_V$ the polynomial in $k[V][X]$ defined as follows:
\begin{equation}
\label{eqn:fVdef}
f_V(X) = \prod_{v \in V} (X - v) \, .
\end{equation}
Recall that Dickson proved there are homogeneous polynomials
$D_s(V)$ for $1 \leq s \leq m$ such that
\begin{equation}
\label{eqn:fV}
f_V(X) = \sum_{s = 0}^m (-1)^s D_s(V) X^{p^{m-s}} \, ,
\end{equation}
where $D_0(V) = 1$. The sequence $D_1(V), \ldots, D_m(V)$ is regular
in $k[V]$, and the invariant ring $k[V]^{\GL(V)}$ is the polynomial algebra
$k[D_1(V), \ldots, D_m(V)]$.

If $\pi \colon V \rightarrow U$ is projection onto a codimension~$\ell$
subspace, then the induced map $k[V] \rightarrow k[U]$ sends
\begin{equation}
\label{eqn:DicksonRes}
D_s(V) \mapsto \begin{cases} D_s(U)^{p^{\ell}} & \text{if $s \leq \dim(U)$,} \\
0 & \text{otherwise.}
\end{cases}
\end{equation}

\section{Existence of polarised systems}
\label{section:existence}
\noindent
For each elementary abelian $p$-group $V$ we shall embed $k[V^*]$ in
$\coho{V}$ by identifying $V^*$ with the image of the Bockstein map
$\coho[1]{V} \rightarrow \coho[2]{V}$.

\subsection{A construction using the norm map}
Let $G$ be a $p$-group of order $p^n$ and $p$-rank $r$ whose centre has
$p$-rank~$z$.
We shall only be interested in the case $r > z$.
Let $U_1, \ldots, U_K$ be representatives of the $G$-orbits in
$\A^C_{z+1}(G)$, which is a $G$-set via the conjugation action.
For each $U \in \A^C_{z+1}(G)$ pick a basis element $x_U$ for the
one-dimensional subspace $\Ann(C)$ of $U^*$, and observe that
$x_U^{p-1}$ is independent of the basis element chosen.
As before, view $U^*$ as a subspace of $\coho[2]{U}$.

Define $\Theta \in \coho{G}$ by
\begin{equation}
\label{eqn:ThetaDef}
\Theta =
\prod_{i = 1}^K \N^G_{U_i}
\left(1 + x_{U_i}^{p-1}\right)^{|G : N_G(U_i)|} \, .
\end{equation}
Now consider the
restriction $\Res_V(\Theta)$ for $V \in \A^C(G)$.
By the Mackey formula
\[
\Res_V (\Theta) =
\prod_{i = 1}^K \prod_{g \in U_i \setminus G / V}
\N^V_{U_i^g \cap V} \, g^* \Res^{U_i}_{U_i \cap {}^g V}
\left(1 + x_{U_i}^{p-1}\right)^{|G : N_G(U_i)|} \, .
\]
The intersection $U_i \cap {}^g V$ always contains $C$,
the largest central elementary abelian subgroup of~$G$. Conversely
the intersection equals $C$ (and $x_{U_i}$ therefore restricts to zero) unless
$U' = U_i^g$ lies in~$V$, in which case
$g^* x_{U_i}^{p-1} = x_{U'}^{p-1}$. Moreover, the number of double cosets
$U_i g V$ leading to this $U'$ is $|N_G(U_i)|/|V|$ and every $U'$~in
$\A^C_{z+1}(V) := \{ U \in \A^C_{z+1}(G) \mid U \leq V \}$ occurs for some~$i$.
So
\begin{equation}
\label{eqn:ResV_Theta}
\Res_V(\Theta) = \prod_{U' \in \A^C_{z+1}(V)}
\N_{U'}^V (1 + x_{U'}^{p-1})^{|G : V|} \, .
\end{equation}
In particular for $U \in \A^C_{z+1}(G)$ one has
\begin{equation}
\label{eqn:ResU_Theta}
\Res_U(\Theta) = 1 + x_U^{(p-1)p^{n-(z+1)}} \, .
\end{equation}
Let $\eta \in \coho[2(p-1)p^{n-(z+1)}]{G}$ be the
homogeneous component of $\Theta$ in this degree.
As the norm map from $\coho{U'}$~to $\coho{V}$ is a ring homomorphism
(see \cite[Proposition~6.3.3]{Evens:book}), we deduce from
Eqn.~\eqref{eqn:ResV_Theta} that
\begin{equation}
\label{eqn:ResV_eta}
\Res_V(\eta) = \hat{\eta}^{|G:V|} \quad \text{for} \quad
\hat{\eta} = \sum_{U' \in \A^C_{z+1}(V)}
\N_{U'}^V x_{U'}^{p-1} \, .
\end{equation}
Denote by $W = W(V)$ the subspace $\Ann(C)$~of $V^*$.
Then $\hat{\eta}$ lies in $k[W]$,
since $\N_{U'}^V (x_{U'})$ is the product of all $\phi \in V^*$
with $\Res_{U'} (\phi) = x_{U'}$.
Moreover $\hat{\eta}$ is by construction a $\GL(W)$-invariant,
so a scalar multiple of $D_1(W)$ for degree reasons. By considering
the restriction to any $U \in \A^C_{z+1}(V)$, we deduce from
Eqn.~\eqref{eqn:DicksonRes} that
\begin{equation}
\label{eqn:ResV_eta_Dickson}
\Res^G_V(\eta) = D_1(W)^{|G : V|} \, .
\end{equation}

\subsection{The existence proof}

\begin{theorem}
\label{theorem:existence}
Let $G$ be a $p$-group of order $p^n$ and $p$-rank $r$ whose centre has
$p$-rank~$z$.
For $1 \leq i \leq z$ define $\zeta_i \in \coho{G}$ by
$\zeta_i = c_{p^n - p^{n-i}}(\rho_G)$, a Chern class of the regular
representation of~$G$.

If $z < r$ define $\eta \in \coho{G}$ as above and 
homogeneous elements $\kappa_1, \kappa_2, \ldots, \kappa_{r-z}
\in \coho{G}$ as follows:
\[
\kappa_j := \rp^{p^{n - z + j - 4}} \cdots \rp^{p^{n - z - 1}}
\rp^{p^{n - z - 2}} (\eta) \in \coho[2(p^{n - z} - p^{n - z - j})]{G}
\]
for $1 \leq j \leq r - z$.
Then for each $1 \leq j \leq r-z$ and for each $V \in \A^C_{z+s}(G)$
one has
\begin{equation}
\label{eqn:existence}
\Res^G_V(\kappa_j) = \begin{cases}
D_j(W)^{|G:V|} & \text{if $j \leq s$, and} \\
0 & \text{otherwise.}
\end{cases}
\end{equation}
Here, $W$~is the subspace of $V^*$ which annihilates~$C$.

Then $\zeta_1, \ldots, \zeta_z, \kappa_1, \ldots, \kappa_{r-z}$ is
a polarised system of parameters for $\coho{G}$. So $\coho{G}$
has both polarised and special polarised systems of parameters.
\end{theorem}

\begin{proof}
Equation~\eqref{eqn:existence} holds for
$\kappa_1 = \eta$ by Eqn.~\eqref{eqn:ResV_eta_Dickson}.
The general case of Eqn.~\eqref{eqn:existence} follows from
Eqn.~\eqref{eqn:DicksonRes} and the action of the Steenrod algebra on the
Dickson invariants (see~\cite{Wilkerson:Dickson}).

Axioms (PS2) and (PS4) follow immediately from Eqn.~\eqref{eqn:existence}.
Axiom (PS3) holds because the Dickson invariants form a regular sequence.
Observe that $\rho_G$ restricts to~$C$ as $p^{n-z}$ copies of the
regular representation $\rho_C$. So the total Chern class $c(\rho_G)$
restricts to $C$ as $c(\rho_C)^{p^{n-z}}$, meaning that
$\Res_C(\zeta_i) = c_{p^z - p^{z-i}}(\rho_C)^{p^{n-z}}$
for $1 \leq i \leq z$.
In view of Eqn.~\eqref{eqn:fVdef} one has $c(\rho_C) = f_{C^*}(1)$
and hence $c_{p^z - p^{z - i}}(\rho_C) = (-1)^i D_i(C^*)$
by Equation~\eqref{eqn:fV}, a well known observation due originally
to Milgram. So $\Res_C(\zeta_i) = (-1)^i D_i(C^*)^{|G:C|}$,
which means that Axiom (PS1) is satisfied and so a polarised
system of parameters has been constructed. Hence special
polarised systems of parameters exist too,
by Lemma~\ref{lemma:polarisedDefinitions}.
\end{proof}

\section{Tightness of Duflot's lower bound}
\label{section:DuflotTight}
\noindent
Recall that Duflot's Theorem~\cite{Duflot:Depth} states that the
depth of $\coho{G}$ is at least~$z$.

\begin{theorem}
\label{theorem:ED1}
Let $G$ be a $p$-group of $p$-rank $r$ whose centre has $p$-rank~$z$.
Then the following statements are equivalent:
\begin{enumerate}
\item
\label{enum:notDetected}
The mod-$p$ cohomology ring $\coho{G}$ is not detected on the family
$\calH^C_{z+1}(G)$.
\item
\label{enum:assocPrime}
$\coho{G}$ has an associated prime $\frakp$ such that the
dimension of $\coho{G}/\frakp$ is $z$.
\item
\label{enum:oneK1}
There is
a polarised system of parameters
$\zeta_1,\ldots,\zeta_z$, $\kappa_1,\ldots,\kappa_{r-z}$
for $\coho{G}$ such that
$\kappa_1$ is a zero divisor in $\coho{G}$.
\item
\label{enum:allK1}
If $\zeta_1,\ldots,\zeta_z$, $\kappa_1,\ldots,\kappa_{r-z}$ is
a polarised system of parameters for $\coho{G}$
then $\kappa_1$ is a zero divisor in $\coho{G}$.
\item
\label{enum:depthz}
The depth of $\coho{G}$ equals $z$.
\end{enumerate}
\end{theorem}

\begin{proof}
Carlson proved in~\cite{Carlson:DepthTransfer} that \eqref{enum:notDetected}
implies \eqref{enum:depthz}\@. A standard commutative algebra argument
shows that \eqref{enum:assocPrime} implies \eqref{enum:depthz}\@.
We saw in Lemma~\ref{lemma:genBrotoHenn} that \eqref{enum:depthz}
implies \eqref{enum:allK1}, and \eqref{enum:oneK1} follows from
\eqref{enum:allK1} by the existence result Theorem~\ref{theorem:existence}.

Now let $\mathfrak{a}$ be a polarised system of parameters
satisfying Statement~\eqref{enum:oneK1}, which is equivalent to $\taua = z$.
So $\tauaH = z$ by Theorem~\ref{theorem:polarisedEqualities}\@.
As in the proof of  that theorem there
is a special polarised system of parameters~$\mathfrak{b}$ which
satisfies $\tauaH[b] = \tauaH = z$, so \eqref{enum:assocPrime}
follows by Lemma~\ref{lemma:kappaAssocPrime}.
Finally consider the definition of~$\tauaH$. If $\tauaH = z$
then $\calH^C_{z+1}(G)$ does not detect $\coho{G}$,
yielding~\eqref{enum:notDetected}\@.
\end{proof}

\section{An example}
Let $G$ be the extraspecial $p$-group $p^{1+2n}_+$ with $n \geq 1$.
This group has order $p^{2n+1}$ and $p$-rank $n+1$. Its centre has
$p$-rank~$1$. If $p$~is odd $G$ has exponent~$p$. The mod-$p$
cohomology ring $\coho{G}$ is Cohen-Macaulay for $p=2$ by Theorem~4.6
of~\cite{Quillen:Extraspecial}, for Quillen computes the cohomology ring
as the quotient of a polynomial algebra by a regular sequence.
Also, Milgram and Tezuka showed in~\cite{MilgramTezuka}
that the cohomology ring is Cohen-Macaulay for $G = 3^{1+2}_+$.

From now on assume that $p$ is odd, with $n \geq 2$ if $p = 3$. Then by
a result of Minh~\cite{Minh:EssExtra}, there are essential classes.
For such groups the centre $C$ is cyclic of order~$p$, and the set
$\calH^C_2(G)$ of centralisers coincides with the set of maximal
subgroups. Consequently $\calH^C_{z+1}(G)$ does not detect $\coho{G}$,
and so $\coho{G}$ has depth~$1$ by the part of Theorem~\ref{theorem:ED1}
proved by Carlson in~\cite{Carlson:DepthTransfer}.

Now let $V$ be a rank $n+1$ elementary abelian subgroup of~$G$,
and let $\hat{\rho}$ be a one-dimensional ordinary representation of~$V$
whose restriction to $C$ is not trivial. Let $\rho$ be the induced
representation of~$G$. Then $\rho$ is an irreducible representation
of degree~$p^n$,
and its character restricts to each $U \in \A^C_{n+1}(G)$ as the sum
of all degree one characters whose restrictions to $C$ coincide with
the restriction of the character of~$\hat{\rho}$. Set
$\zeta_1 := c_{p^n}(\rho)$ and $\kappa_j := c_{p^n - p^{n-j}}(\rho)$
for $1 \leq j \leq n$. Then $\zeta_1, \kappa_1, \ldots, \kappa_n$
satisfies the axioms for a polarised system of parameters for $\coho{G}$.
So combining Minh's result with Theorem~\ref{theorem:ED1}
one deduces that $\kappa_1$ has nontrivial annihilator in $\coho{G}$.
Conversely a direct proof of this fact would yield a new proof
of Minh's result. If $n=1$ and $p > 3$ it is known (see~\cite{Leary:integral})
that $c_2(\rho)$ is a nonzero essential class which annihilates~$\kappa_1$.


\begin{thebibliography}{10}

\bibitem{Benson:ImTr}
D.~J. Benson.
\newblock The image of the transfer map.
\newblock {\em Archiv Math.}, 61:7--11, 1993.

\bibitem{Benson:PolyInvts}
D.~J. Benson.
\newblock {\em Polynomial invariants of finite groups}.
\newblock London Math. Soc. Lecture Note Series, vol.~190. Cambridge University
  Press, Cambridge, 1993.

\bibitem{BrHe:Central}
C.~Broto and H.-W. Henn.
\newblock Some remarks on central elementary abelian $p$-subgroups and
  cohomology of classifying spaces.
\newblock {\em Quart. J. Math. Oxford Ser. (2)}, 44(174):155--163, 1993.

\bibitem{Carlson:DepthTransfer}
J.~F. Carlson.
\newblock Depth and transfer maps in the cohomology of groups.
\newblock {\em Math. Z.}, 218(3):461--468, 1995.

\bibitem{Carlson:Problems}
J.~F. Carlson.
\newblock Problems in the calculation of group cohomology.
\newblock In P.~Dr{\"a}xler, G.~O. Michler, and C.~M. Ringel, editors, {\em
  Computational methods for representations of groups and algebras (Essen,
  1997)}, pages 107--120. Birkh\"auser, Basel, 1999.

\bibitem{Duflot:Depth}
J.~Duflot.
\newblock Depth and equivariant cohomology.
\newblock {\em Comment. Math. Helv.}, 56(4):627--637, 1981.

\bibitem{Evens:book}
L.~Evens.
\newblock {\em The cohomology of groups}.
\newblock Oxford Univ.\@ Press, Oxford, 1991.

\bibitem{Leary:integral}
I.~J. Leary.
\newblock The integral cohomology rings of some $p$-groups.
\newblock {\em Math. Proc. Cambridge Philos. Soc.}, 110(1):25--32, 1991.

\bibitem{MilgramTezuka}
R.~J. Milgram and M.~Tezuka.
\newblock The geometry and cohomology of ${M}\sb {12}$. {I}{I}.
\newblock {\em Bol. Soc. Mat. Mexicana (3)}, 1(2):91--108, 1995.

\bibitem{Minh:EssExtra}
P.~A. Minh.
\newblock Essential cohomology and extraspecial $p$-groups.
\newblock {\em Trans. Amer. Math. Soc.}, 353(5):1937--1957, 2000.

\bibitem{Quillen:Extraspecial}
D.~Quillen.
\newblock The mod-2 cohomology rings of extra-special 2-groups and the spinor
  groups.
\newblock {\em Math. Ann.}, 194:197--212, 1971.

\bibitem{Wilkerson:Dickson}
C.~Wilkerson.
\newblock A primer on the {D}ickson invariants.
\newblock In H.~R. Miller and S.~B. Priddy, editors, {\em Proceedings of the
  Northwestern Homotopy Theory Conference (Evanston, Ill., 1982)}, Contemporary
  Math., vol.~19, pages 421--434, Providence, RI, 1983. Amer. Math. Soc.

\end{thebibliography}
\end{document}